\documentstyle[12pt]{article}

\makeatletter

\def\newremark#1{\@ifnextchar[{\@ormkwsq{#1}}{\@nrmkwsq{#1}}}

\def\@nrmkwsq#1#2{%
\@ifnextchar[{\@xnrmkwsq{#1}{#2}}{\@ynrmkwsq{#1}{#2}}}

\def\@xnrmkwsq#1#2[#3]{\expandafter\@ifdefinable\csname #1\endcsname
{\@definecounter{#1}\@addtoreset{#1}{#3}%
\expandafter\xdef\csname the#1\endcsname{\expandafter\noexpand
  \csname the#3\endcsname \@rmkcountersep \@rmkcounter{#1}}%
\global\@namedef{#1}{\@rmkwsq{#1}{#2}}
\global\@namedef{end#1}{\@endremarkwithsquare}}}

\def\@ynrmkwsq#1#2{\expandafter\@ifdefinable\csname #1\endcsname
{\@definecounter{#1}%
\expandafter\xdef\csname the#1\endcsname{\@rmkcounter{#1}}%
\global\@namedef{#1}{\@rmkwsq{#1}{#2}}
\global\@namedef{end#1}{\@endremarkwithsquare}}}

\def\@ormkwsq#1[#2]#3{\expandafter\@ifdefinable\csname #1\endcsname
  {\global\@namedef{the#1}{\@nameuse{the#2}}%
\global\@namedef{#1}{\@rmkwsq{#2}{#3}}%
\global\@namedef{end#1}{\@endremarkwithsquare}}}

\def\@rmkwsq#1#2{\refstepcounter
    {#1}\@ifnextchar[{\@yrmkwsq{#1}{#2}}{\@xrmkwsq{#1}{#2}}}

\def\@xrmkwsq#1#2{\@beginremark{#2}{\csname the#1\endcsname}\ignorespaces}
\def\@yrmkwsq#1#2[#3]{\@opargbeginremark{#2}{\csname
       the#1\endcsname}{#3}\ignorespaces}

\def\@rmkcounter#1{\noexpand\arabic{#1}}
\def\@rmkcountersep{.}
\def\@beginremark#1#2{\trivlist \item[\hskip \labelsep{\bf #1\ #2.}]}
\def\@opargbeginremark#1#2#3{\trivlist
      \item[\hskip \labelsep{\bf #1\ #2\ (#3)}]}
\def\@endremarkwithsquare{~\hspace{\fill}~$\square$\endtrivlist}

\def\@eqnnum{\hbox to .01pt{}\rlap{\rm \hskip -\displaywidth(\theequation)}}

   \def\@begintheorem#1#2{\sl \trivlist \item[\hskip \labelsep{\bf #2\ #1}]}
   \def\@opargbegintheorem#1#2#3{\sl \trivlist
            \item[\hskip \labelsep{\bf #2\ #1\ (#3)}]}

   \def\section{\@startsection {section}{1}{\z@}{-3.5ex plus -1ex minus
    -.2ex}{2.3ex plus .2ex}{\large\bf}}

   \def\subsection{\@startsection{subsection}{2}{\z@}{-3.25ex plus -1ex minus
   -.2ex}{1.5ex plus .2ex}{\normalsize\bf}}

\newenvironment{eqn}{\refstepcounter{subsection}
$$}{\leqno{\rm(\thesubsection)}$$\global\@ignoretrue}

\newenvironment{subeqn}{\refstepcounter{subsubsection}
$$}{\leqno{\rm(\thesubsubsection)}$$\global\@ignoretrue}

\def\@rmkcounter#1{\noexpand\arabic{#1}}
\def\@rmkcountersep{.}
\def\@beginremark#1#2{\trivlist \item[\hskip \labelsep{\bf #2\ #1.}]}
\def\@opargbeginremark#1#2#3{\trivlist
      \item[\hskip \labelsep{\bf #2\ #1\ (#3).}]}
\def\@endremarkwithsquare{~\hspace{\fill}~$\Box$\endtrivlist}
\makeatother

\newenvironment{prf}[1]{\trivlist
\item[\hskip \labelsep{\it
#1.\hspace*{.3em}}]}{~\hspace{\fill}~$\Box$\endtrivlist}

\newenvironment{proof}{\begin{prf}{\bf Proof}}{\end{prf}}

\let\tempcirc=\circ
\def\circ{\mathord{\raise0.25ex\hbox{$\scriptscriptstyle\tempcirc$}}}

\newcommand{\ZZ}{{\bf Z}}
\newcommand{\QQ}{{\bf Q}}
\newcommand{\RR}{{\bf R}}
\newcommand{\CC}{{\bf C}}
\newcommand{\FF}{{\bf F}}
\newcommand{\HH}{{\bf H}}
\newcommand{\PP}{{\bf P}}
\newcommand{\SL}{{\rm SL}}
\newcommand{\PGL}{{\rm PGL}}
\newcommand{\GL}{{\rm GL}}
\newcommand{\Qbar}{\overline{\QQ}}
\newcommand{\End}{{\rm End}}
\newcommand{\Lie}{{\rm Lie}}
\newcommand{\Gal}{{\rm Gal}}
\newcommand{\Pic}{{\rm Pic}}
\newcommand{\discr}{{\rm discr}}
\newcommand{\Spec}{{\rm Spec}}
\newcommand{\rH}{{\rm H}}

\newcommand{\eps}{\varepsilon}
\newcommand{\Li}{{\rm Li}}
\newcommand{\pr}{{\rm pr}}
\newcommand{\SO}{{\rm SO}}
\newcommand{\Hbar}{\overline{H}}

\newtheorem{theorem}[subsection]{Theorem.}
\newtheorem{proposition}[subsection]{Proposition.}
\newtheorem{lemma}[subsection]{Lemma.}

\newremark{remark}[subsection]{Remark}
\newremark{subremark}[subsubsection]{Remark}
\newremark{notation}[subsection]{Notation}
\newremark{example}[subsection]{Example}
\newremark{definition}[subsection]{Definition}

\begin{document}
\title{Special points on the product of two modular curves.}
\author{Bas Edixhoven\footnote{partially supported by the
Institut Universitaire de France}}
\maketitle

\section{Introduction.}\label{section1}
It is well known that the $j$-invariant establishes a bijection
between $\CC$ and the set of isomorphism classes of elliptic curves
over~$\CC$, see for example \cite{Silverman1}. The endomorphism ring of an
elliptic curve $E$ over $\CC$ is either $\ZZ$ or an order in an imaginary
quadratic extension of $\QQ$; in the second case $E$ is said to be a CM
elliptic curve (CM meaning complex multiplication).
A complex number $x$ is said to be CM if the corresponding elliptic curve
over $\CC$ is~CM. A point $(x_1,x_2)$ in $\CC^2$ is defined to be CM if both
$x_1$ and $x_2$ are~CM. The aim of this article is to determine all
irreducible
algebraic curves $C$ in $\CC^2$ containing infinitely many CM points. In other
words, we want to determine all irreducible polynomials $f$ in $\CC[x_1,x_2]$
that vanish at infinitely many CM points. The motivation for doing this comes
from a conjecture of Frans Oort (see \cite[Chapter~IV, \S1]{Moonen1} for a
precise statement), saying roughly that the irreducible components of the
Zariski closure of any set of CM points in any Shimura variety are sub
Shimura varieties. For the irreducible components of dimension zero this
is trivially true. For those of dimension one Oort's conjecture was in fact
stated earlier by Yves Andr\'e as a problem in \cite[Chapter~X, \S1]{Andre2}.

We view $\CC^2$ as the Shimura variety which is the moduli space of pairs of
elliptic curves. Then the irreducible sub Shimura varieties of dimension one
are the following: $\CC\times\{x_2\}$ with $x_2$ a CM point,
$\{x_1\}\times\CC$ with $x_1$ a CM point, or the image in $\CC^2$, under the
usual map, of the modular curve $Y_0(n)$ for some integer $n\geq1$. Recall
that, for $n\geq1$, $Y_0(n)$ is the modular curve classifying elliptic curves
with a cyclic subgroup of order $n$, or, equivalently, cyclic isogenies of
degree $n$ between elliptic curves. The usual map from $Y_0(n)$ to $\CC^2$
sends an isogeny to its source and target, i.e., $\phi\colon E_1\to E_2$
is sent to $(j(E_1),j(E_2))$. We will prove the following result, giving
evidence for the conjecture just mentioned.

\begin{theorem}\label{thm1.1}
Assume the generalized Riemann hypothesis for imaginary quadratic fields.
Let $C$ be an irreducible algebraic curve in\/ $\CC^2$ containing infinitely
many CM points and such that neither of its projections to\/ $\CC$ is
constant. Then $C$ is the image of\/ $Y_0(n)$ for some $n\geq1$.
\end{theorem}
\begin{remark}\label{rmk1.2}
In the proof of Theorem~\ref{thm1.1} we will see that the state of the art
in analytic number theory is such that the Riemann hypothesis is ``almost
not needed'' (see Remark~\ref{rmk5.4}).
It is clear that Theorem~\ref{thm1.1} implies similar statements for
curves contained in the product of two modular curves. In particular,
if one assumes GRH, Oort's conjecture is true for curves contained in the
product of two modular curves.
\end{remark}
\begin{remark}\label{rmk1.3}
Ben Moonen has proved Oort's conjecture for the sets of CM points in moduli
spaces of abelian varieties such that there exists a prime number $p$ at which
all the CM points are canonical in the sense that they have an ordinary
reduction of which they are the Serre-Tate canonical lift (see
\cite[Chapter~IV, \S1]{Moonen1}). Yves Andr\'e has proved the conclusion of
Theorem~\ref{thm1.1} with the Riemann hypothesis replaced by the assumption
that the Zariski closure of $C$ in $\PP^1\times\PP^1$ meets
$\{\infty\}\times\CC$ only in points $(\infty,x_2)$ with $x_2$ a CM point
(see \cite{Andre1}). In the case where $C$ meets the union of
$\{\infty\}\times\CC$ and $\CC\times\{\infty\}$ only in $\infty\times\infty$
he has a very simple proof.
\end{remark}
The idea of the proof of Theorem~\ref{thm1.1} is the following. We use the
Galois action on the set of CM $j$-invariants to show that for all but
finitely many CM points $(x_1,x_2)$ on $C$ the CM fields of $x_1$ and $x_2$
coincide. Then we consider intersections of $C$ with its images under
certain Hecke operators. The Riemann hypothesis implies that $C$ is actually
contained in some of these images. To finish, we consider an irreducible
component $X$ of the inverse image of $C$ in $\HH\times\HH$, the product of
the complex upper half plane by itself, and show that the stabilizer of
$X$ in $\SL_2(\RR)\times\SL_2(\RR)$ is of the kind it should be.

\section{Some facts about CM elliptic curves.}\label{section2}
Before we start with the proof of Theorem~\ref{thm1.1}, we need to
recall some facts about CM elliptic curves. These facts can be found
for example in \cite[Appendix~C, \S11]{Silverman1}.
First of all, CM elliptic curves are defined over~$\Qbar$.
Let $K$ be an imaginary quadratic extension of $\QQ$, with a given
embedding in~$\Qbar$. Let $O_K\subset K$ be the ring of integers.
Every subring $A$ of $O_K$ of finite index is of the form $O_{K,f}:=\ZZ+fO_K$
for a unique integer $f\geq1$. For $f\geq1$ let $S_{K,f}$ be the set of
isomorphism classes of pairs $(E,\alpha)$, with $E$ an elliptic curve over
$\Qbar$ and $\alpha\colon O_{K,f}\to\End(E)$ an isomorphism of rings inducing
the given embedding of $K$ into $\Qbar$ via the action on $\Lie(E)$. The
group $G_K:=\Gal(\Qbar/K)$ acts on $S_{K,f}$. But also the Picard group
$\Pic(O_{K,f})$ acts on $S_{K,f}$ by the following formula:
\begin{eqn}\label{eqn2.1}
(E,[L]) \mapsto E\otimes_{O_{K,f}}L,
\end{eqn}
where $L$ is an invertible $O_{K,f}$-module, $[L]$ its equivalence class
and $E\otimes_{O_{K,f}}L$ the cokernel of the map $p\colon E^2\to E^2$ if
$p\colon O_{K,f}^2\to O_{K,f}^2$ has cokernel $L$ (view $p$ as a matrix with
coefficients in $O_{K,f}$). If we choose an embedding
of $\Qbar$ in $\CC$ and write $E(\CC)$ as $\CC$ modulo a lattice $\Lambda$,
then $(E\otimes_{O_{K,f}}L)(\CC)$ is the quotient of $\CC\otimes_{O_{K,f}}L$
by $\Lambda\otimes_{O_{K,f}}L$. The actions by $G_K$ and $\Pic(O_{K,f})$ on
$S_{K,f}$ commute.
\begin{proposition}\label{prop2.2}
The set $S_{K,f}$ is a $\Pic(O_{K,f})$-torsor, i.e., the action of\/
$\Pic(O_{K,f})$ is free and has exactly one orbit.
\end{proposition}
\begin{proof} (Sketch.) For every $(E,\alpha)$ and $\Lambda$ as above,
$\End_{O_{K,f}}(\Lambda)=O_{K,f}$. Moreover, $O_{K,f}$ is of the form
$\ZZ[x]/(g)$. It follows that $\Lambda$ is an invertible $O_{K,f}$-module.
\end{proof}
It follows that $G_K$ acts on $S_{K,f}$ via a morphism $G_K\to\Pic(O_{K,f})$.
This morphism is surjective and unramified outside~$f$. The Frobenius
element at a maximal ideal $m$ not containing $f$ is the element $[m]^{-1}$
of $\Pic(O_{K,f})$ (all this can be seen from deformation theory, using the
theorem of Serre-Tate, or from class field theory). Let $H_{K,f}$ be the
Galois extension of $K$ corresponding to this quotient $\Pic(O_{K,f})$
of~$G_K$. We remark that we have $H_{K,f}=K(j(E))$ for all $(E,\alpha)$
in~$S_{K,f}$.

\section{The two CM fields are almost always equal.}\label{section3}
Let $C_\CC\subset\CC^2$ be as in Theorem~\ref{thm1.1} (i.e., it is
irreducible,
it contains infinitely many CM points and its two projections to $\CC$
are not constant). Since all CM points have coordinates in $\Qbar$, $C_\CC$ is
defined over $\Qbar$, in the sense that it is the locus of  zeros of an
irreducible polynomial, call it $f$, with coefficients in~$\Qbar$.
It will be convenient for us to work with a curve defined over $\QQ$, hence
we let $C$ be the union of the finitely many conjugates of~$C_\CC$. Then
$C$ is defined by the product $F$ of the Galois conjugates of $f$, if we take
$f$ such that it has a non-zero coefficient in~$\QQ$. Let $d_1$ and $d_2$ be
the degrees of $F$ with respect to the second and first variable.
Then $d_i$ is the degree of the $i$th projection from $C$ to~$\CC$.
For $x$ in $\CC$ we will denote the endomorphism ring of the corresponding
elliptic curve by~$\End(x)$. For a CM point $x$ in $\CC$ we will call
$\QQ\otimes\End(x)$ the CM field of~$x$. Note that the isogeny class of a
CM elliptic curve over $\Qbar$ consists of all elliptic curves with the
same CM field.
We want to prove that $C$ is the image in $\CC^2$ of some~$Y_0(n)$. Our
first step in this direction is the following proposition.

\begin{proposition}\label{prop3.1}
Let $C$ be as above. For all but finitely many CM points $(x_1,x_2)$ in $C$
the CM fields of $x_1$ and $x_2$ coincide.
\end{proposition}
\begin{proof}
Suppose that $(x_1,x_2)$ is a CM point in $C(\Qbar)$ such that the two
CM fields $K_1$ and $K_2$ are different. Since $C$ is defined over $\QQ$,
$\QQ(x_1,x_2)$ has degree at most $d_2$ over $\QQ(x_1)$ and degree at
most $d_1$ over~$\QQ(x_2)$. Let $L$ be the field generated by $K_1$
and $K_2$, and $M$ the intersection of $L(x_1)$ and~$L(x_2)$. Let us write
$\End(x_i)=O_{K_i,f_i}$ for $i=1$ and~$2$. The field $L(x_i)$ is an abelian
Galois extension of $L$, of degree at least $|\Pic(O_{K_i,f_i})|/2$.
The degrees of $L(x_1,x_2)$ over $L(x_2)$ and $L(x_1)$ are equal to those of
$L(x_1)$ and $L(x_2)$ over $M$, respectively. This gives us:
\begin{eqn}\label{eqn3.2}
|\Pic(O_{K_i,f_i})| \leq 2d_i[M:L].
\end{eqn}
We will now work to get a suitable upper bound for~$[M:L]$. The group
$\Gal(L(x_1,x_2)/\QQ)$ is an extension of $\Gal(L/\QQ)$ by the abelian group
$\Gal(L(x_1,x_2)/L)$. Hence the action of $\Gal(L(x_1,x_2)/\QQ)$ on
$\Gal(L(x_1,x_2)/L)$
by conjugation factors through an action of $\Gal(L/\QQ)$. In the same way,
$\Gal(L/\QQ)$ acts on the two groups $\Gal(L(x_i)/L)$, which we view
as subgroups of $\Gal(K_i(x_i)/K_i)$. Now $\Gal(L/\QQ)$
is equal to $\Gal(K_1/\QQ)\times\Gal(K_2/\QQ)$, hence equal to
$\ZZ/2\ZZ\times\ZZ/2\ZZ$. The action of $\Gal(L/\QQ)$ on $\Gal(L(x_i)/L)$
factors through $\Gal(K_i/\QQ)$ and as such coincides with the restriction of
the action of $\Gal(K_i/\QQ)$ on $\Gal(K_i(x_i)/K_i)=\Pic(O_{K_i,f_i})$.
\begin{lemma}\label{lemma3.3}
Let $K$ be a quadratic imaginary field and $f\geq1$. Then the non-trivial
element $\sigma$ of\/ $\Gal(K/\QQ)$ acts as $-1$ on\/ $\Pic(O_{K,f})$.
\end{lemma}
\begin{proof}
The endomorphism $\sigma+1$ of $\Pic(O_{K,f})$ factors through the norm
map from $\Pic(O_{K,f})$ to~$\Pic(\ZZ)$.
\end{proof}
Now note that $\Gal(M/L)$ is a quotient of both $\Gal(L(x_i)/L)$, so
the action of $\Gal(L/\QQ)$ on it is by the non-trivial character given by
the first projection, but also by the second projection. This implies that
$\Gal(M/L)$ is killed by multiplication by two.

\begin{lemma}\label{lemma3.4}
Let $K$ be an imaginary quadratic field and $f\geq1$. Then the dimension
of the $\FF_2$-vector space $\Pic(O_{K,f})\otimes\FF_2$ is at most the number
of odd primes dividing the discriminant $\discr(O_{K,f})$ of $O_{K,f}$ plus
ten.
\end{lemma}
\begin{proof}
(Sketch.) The exact bound we give does not matter so much, so we just give
some indications. First one notes that there is an exact sequence:
\begin{subeqn}\label{eqn3.4.1}
(K\otimes\QQ_2)^* \to \Pic(O_{K,f}) \to \Pic(O_{K,f}[1/2]) \to 0 .
\end{subeqn}
Let $S:=\Spec(O_{K,f}[1/2])$ and $T:=\Spec(\ZZ[1/2])$.
The Kummer sequence gives a surjection from $\rH^1(S_{\rm et},\FF_2)$
onto the $2$-torsion subgroup of $\Pic(S)$, which has the same dimension
as $\Pic(S)\otimes\FF_2$. One deals with $\rH^1(S_{\rm et},\FF_2)$ by
projecting to~$T_{\rm et}$.
\end{proof}
Since $\Gal(M/L)$ is killed by 2 and a quotient of a subgroup of
$\Pic(O_{K_i,f_i})$, we have:
\begin{eqn}\label{eqn3.5}
\log_2 [M:L] \leq |\{2\neq p | \discr(O_{K_i,f_i})\}| + 10,
\qquad i\in\{1,2\}.
\end{eqn}
On the other hand, we have Siegel's theorem (see \cite{Oesterle1}), stating
that:
\begin{eqn}\label{eqn3.6}
\log |\Pic(O_{K_i,f_i})| = (1/2 + {\rm o}(1))\log |\discr(O_{K_i,f_i})|,
\qquad (|\discr(O_{K_i,f_i})|\to\infty).
\end{eqn}
Combining equations (\ref{eqn3.5}) and (\ref{eqn3.6}) shows that
$|\Pic(O_{K_i,f_i})|/[M:L]$ tends to infinity as the discrminiant of
$O_{K_i,f_i}$ tends to infinity. But then equation (\ref{eqn3.2}) can hold
for only finitely many $(x_1,x_2)$. This ends the proof of
Proposition~\ref{prop3.1}.
\end{proof}

\begin{remark}\label{rmk3.7}
The proof of Proposition~\ref{prop3.1} shows actually more: the function
on the set of CM points on $C$ that sends $(x_1,x_2)$ to $f_1/f_2$ takes
only finitely many values. Using this, one can reduce the proof of
Theorem~\ref{thm1.1} to the case where there are infinitely many CM points
$(x_1,x_2)$ on $C$ with $\End(x_1)=\End(x_2)$ (one replaces $C$ by its
image under a suitable Hecke correspondence). As we do not know how to
exploit this, we do not go into further detail.
\end{remark}
\begin{remark}\label{rmk3.8}
Proposition~\ref{prop3.1} was also proved by Yves Andr\'e in \cite{Andre1},
and also by Ching-Li Chai (not published).
\end{remark}

\section{Intersecting $C$ with something.}
\label{section4}
We continue the proof of Theorem~\ref{thm1.1}. So we let $C$ be
as before. At this point we already know that we have infinitely many
CM points $(x_1,x_2)$ on $C$ for which $x_1$ and $x_2$ are isogeneous
because they have the same CM field. We have to prove that there
is an integer $n\geq1$ such that for infinitely many $(x_1,x_2)$ there
exists an isogeny of degree $n$ between $x_1$ and~$x_2$. A direct
approach for this is the following. Consider a CM point $(x_1,x_2)$ such that
$x_1$ and $x_2$ have the same CM field, say $K$, and an isogeny from
$x_1$ to $x_2$ of minimal degree, say~$n$. One can get an upper bound
for $n$ in terms of the discriminants of the~$\End(x_i)$.
By Remark~\ref{rmk3.7}, one can assume that $\End(x_1)=\End(x_2)=O_{K,f}$
and get an upper bound for $n$ from Minkowski's theorem on ideals of small
norm representing elements of the class group; the bound is a constant times
$|\discr(O_{K,f})|^{1/2}$. Then one considers the intersection of $C$
with~$Y_0(n)$. The degrees of both projections from $Y_0(n)$ to $\CC$
are equal to $\psi(n)$, where $\psi(n)=n\prod_{p|n}(1+1/p)$.
The Picard group of $\PP^1\times\PP^1$ (over a field, say $\QQ$) is
isomorphic to $\ZZ\times\ZZ$, the isomorphism sending an effective divisor to
the degrees of its two projections to~$\PP^1$. The intersection form is
the following: $(a,b)\cdot(c,d)=ad+bc$. Hence the intersection number
of the Zariski closures in $\PP^1\times\PP^1$ of $C$ and $Y_0(n)$
is $\psi(n)(d_1+d_2)$. Since both curves we intersect are defined over
$\QQ$, the intersection contains all Galois conjugates of $(x_1,x_2)$, of
which there are~$|\Pic(O_{K,f})|$.
So if $|\Pic(O_{K,f})|$ exceeds $\psi(n)(d_1+d_2)$, the proof is finished,
since then the intersection is not proper. Unfortunately,
equation (\ref{eqn3.6}) does not imply such an inequality.

Nevertheless, the idea of intersecting $C$ with something is a good
one. Natural ``somethings'' to take are images of $C$ itself under Hecke
correspondences. Again, we consider a CM point $(x_1,x_2)$ on $C$
such that the CM fields of $x_1$ and $x_2$ coincide. Let $K$, $f_1$ and
$f_2$ be defined by: $\End(x_i)=O_{K,f_i}$. Let $f$ be the least common
multiple of $f_1$ and~$f_2$. One easily checks that the field generated by
$H_{K,f_1}$ and
$H_{K,f_2}$ is $H_{K,f}$, hence the orbit of $(x_1,x_2)$ under the
action of $G_K$ is a $\Gal(H_{K,f}/K)$-torsor. Recall from \S\ref{section2}
that we can identify $\Gal(H_{K,f}/K)$ with~$\Pic(O_{K,f})$.
For $\sigma$ in $\Gal(H_{K,f}/K)$ corresponding to the class $[I]$ of
an invertible ideal $I$ of $O_{K,f}$, there are isogenies from $x_1$ to
$\sigma(x_1)$ and from $x_2$ to $\sigma(x_2)$ whose kernels are
isomorphic, as $O_{K,f}$-modules, to~$O_{K,f}/I$. Hence if we take $I$
such that $O_{K,f}/I$ is a cyclic group of some order $n$, then
$\sigma(x_i)$ is in $T_n(x_i)$ for $i$ equals 1 and~2, where $T_n$ is
the correspondence on $\CC$ that sends an elliptic curve to the sum (as
divisors) of its quotients by its cyclic subgroups of order~$n$.
(Let us note that this $T_n$ is not the same as the correspondence on $\CC$
that is usually called $T_n$ if $n$ is not square free, since the usual one
involves a sum over all subgroups of order~$n$.)
Let $T_n\times T_n$ be the correspondence on $\CC\times\CC$ that is the
product of $T_n$ on each factor: it sends a pair $(E_1,E_2)$ of
elliptic curves to the sum of the $(E_1/G_1,E_2/G_2)$, where $G_i$ is a
cyclic subgroup of order $n$ in~$E_i$. Then $(x_1,x_2)$ is in the
intersection of $C$ and $(T_n\times T_n)C$, because $x_i$ is in
$T_n(\sigma(x_i))$ and $(\sigma(x_1),\sigma(x_2))$ is in~$C$. Since
both $C$ and $(T_n\times T_n)C$ are defined over $\QQ$, their intersection
contains all Galois conjugates of~$(x_1,x_2)$. Hence the intersection
has at least $|\Pic(O_{K,f})|$ elements. Let us now calculate the degrees
of the projections of $(T_n\times T_n)C$ to~$\CC$. By definition,
$(T_n\times T_n)C$ consists of the $(x,y)$ such that there
exist $u$ and $v$ in $\CC$ with $(u,v)$ in $C$, and cyclic isogenies of
degree $n$ from $u$ to $x$ and from $v$ to~$y$. Let $x$ be in~$\CC$.
Then there are $\psi(n)$ $u$'s with $x\in T_n(u)$. For each such a $u$
there are $d_1$ $v$'s with $(u,v)$ on $C$. For each such a $v$ there are
$\psi(n)$ $y$'s in~$T_n(v)$. This shows that the degree of the first
projection of $(T_n\times T_n)C$ is $\psi(n)^2d_1$. Of course, for the
second projection one has the analogous result. So, for the intersection
number of $C$ and $(T_n\times T_n)C$ we find $2d_1d_2\psi(n)^2$.
We conclude that if $|\Pic(O_{K,f})|$ is bigger than $2d_1d_2\psi(n)^2$,
then $C$ is contained in $(T_n\times T_n)C$. The next thing to do is to
see if there do exist ideals $I$ with the required properties.

Let $x_1$, $x_2$, $K$ and $f$ be as above. Let $p$ be a prime number
that splits in $O_{K,f}$, i.e., such that $O_{K,f}\otimes\FF_p$ is
isomorphic to $\FF_p\times\FF_p$. For $I$ we take one of the two maximal
ideals containing~$p$. As explained above, we have the following
implication:
\begin{eqn}\label{eqn4.1}
2d_1d_2(p+1)^2 < |\Pic(O_{K,f})| \quad\mbox{\rm implies} \quad
C\subset (T_p\times T_p)C .
\end{eqn}
Equation~(\ref{eqn3.6}) tells us that
$|\Pic(O_{K,f})|=|\discr(O_{K,f})|^{1/2+{\rm o}(1)}$. So we want $p$ to
be at most something as $|\discr(O_{K,f})|^{1/4}$. More precisely:

\begin{proposition}\label{prop4.2}
Suppose that there exists $\eps>0$ such that, when $K$ ranges through all
imaginary quadratic fields and $f$ through all positive integers, the
number of primes $p<|\discr(O_{K,f})|^{1/4-\eps}$
that are split in $O_{K,f}$ tends to infinity as $|\discr(O_{K,f})|$ tends to
infinity. Then there are infinitely many primes $p$ such that $C$ is
contained in $(T_p\times T_p)C$.
\end{proposition}
\begin{proof}
Because we have infinitely many CM points $(x_1,x_2)$ on $C$, we know that
the discriminants $|\discr(O_{K,f})|$ associated to them as above tend to
infinity. The implication (\ref{eqn4.1}) and equation (\ref{eqn3.6}) give
us the infinitely many required primes.
\end{proof}

\section{Existence of small split primes.} \label{section5}
The aim of this section is to prove the hypothesis in
Proposition~\ref{prop4.2}. It turns out that this is no problem at all
if one assumes GRH for imaginary quadratic fields and uses the resulting
effective Chebotarev theorem of Lagarias, Montgomery and Odlyzko as
stated in~\cite{Serre1}.

For $K$ an imaginary quadratic field, $f$ a positive integer and $x\geq2$ a
real number, let $\pi_{K,f}(x)$ be the number of primes $p\leq x$ that
are split in~$O_{K,f}$, let $d_K:=|\discr(O_K)|$
and let $d_{K,f}:=|\discr(O_{K,f})|$. Note that $d_{K,f}=f^2d_K$. As usual,
let $\Li(x):=\int_2^xdt/\log(t)$.
Theorem~4 of \cite{Serre1} and the second remark following it say that,
for $x$ sufficiently big and for all $K$ as above for which GRH holds,
one has:
\begin{eqn}\label{eqn5.1}
\left|\pi_{K,1}(x)-{1\over2}\Li(x)\right| \leq
{1\over6}x^{1/2}\left(\log(d_K) + 2\log(x)\right) .
\end{eqn}
Since the number of primes dividing $f$ is at most $\log_2(f)$, equation
(\ref{eqn5.1}) implies:
\begin{eqn}\label{eqn5.2}
\pi_{K,f}(x) \geq {x\over2\log(x)}
\left(\Li(x){\log(x)\over x} - {\log(x)\over 3x^{1/2}}
\left(\log(d_K)+2\log(x)\right) - {2\log(x)\log(f)\over x\log(2)}\right).
\end{eqn}
If $x$ tends to infinity, $\Li(x)\log(x)/x$ tends to 1 and
$\log(x)^2/x^{1/2}$ tends to~0. One checks easily that for $x$ sufficiently
big (i.e., bigger than some absolute constant), and bigger than
$\log(d_{K,f})^2(\log(\log(d_{K,f}))^2$, one has
$\log(x)\log(d_K)/3x^{1/2}<c<1$, with $c$ independent of $K$ and~$f$.
Under the same conditions, $\log(x)\log(f)/x$ tends to zero if $x$ tends to
infinity.
This means that we have proved the following proposition.
\begin{proposition}\label{prop5.3}
Let $C$ be as before (i.e., as in the beginning of~\S\ref{section3}).
Assume GRH for all imaginary quadratic fields.
Then there exist infinitely many primes $p$ such that $C$ is contained
in $(T_p\times T_p)C$.~\hfill$\Box$
\end{proposition}
\begin{remark}\label{rmk5.4}
Of course, the question remains whether one can prove the hypothesis of
Proposition~\ref{prop4.2} without assuming GRH. Etienne Fouvry tells me
the following. He shows that for $r>0$ and all $n$, the set of
$d_{K,f}$ such that the number of primes $p<d_{K,f}^r$ that are split
in $O_{K,f}$ is at most $n$, has density zero (i.e., the number of such
$d_{K,f}<x$ is ${\rm o}(x)$ for $x\to\infty$). Moreover, he says that the
exponent $1/4$ is critical, in the sense that one can prove that for
all $\eps>0$, the number of primes $p<d_{K,f}^{1/4+\eps}$ that are split
in $O_{K,f}$ tends to infinity as $d_{K,f}$ tends to infinity. To prove
this, he uses a result of Linnik and Vinogradov in  \cite{LinnikVinogradov},
see also \cite{Friedlander}. The central point in \cite{LinnikVinogradov}
is an upper bound for short character sums by Burgess, in which the exponent
$1/4+\eps$ appears. This $1/4$ has not moved in the last 30 years.
\end{remark}

\section{Some topological arguments.}\label{section6}
In this section we finish the proof of Theorem~\ref{thm1.1} by combining
Proposition~\ref{prop5.3} with the following theorem, which gives yet
another characterization of modular curves.

\begin{theorem}\label{thm6.1}
Let $C$ in $\CC^2$ be an irreducible algebraic curve. Let $d_1$ and $d_2$ be
the degrees of its two projections to~$\CC$. Suppose that $d_1$ and $d_2$
are both non-zero, and that we have $C\subset (T_n\times T_n)C$ for some
square free integer $n>1$ that is composed of primes $p\geq \max\{5,d_1\}$.
Then $C$ is the image of\/ $Y_0(m)$ in\/ $\CC^2$ for some $m\geq1$.
\end{theorem}
Let us first show that this theorem and Proposition~\ref{prop5.3}
imply Theorem~\ref{thm1.1}. So let $C_\CC$ and $C$ be as in the beginning
of~\S\ref{section3}. Recall that $C$ is the union of the finitely many
Galois conjugates of the irreducible component $C_\CC$ of it.
We know that there are infinitely many primes $p$ such that $C$ is contained
in $(T_p\times T_p)C$. For such a prime $p$, let $T_{C,p}$ denote the
correspondence on $C$ induced by $T_p\times T_p$. By this we mean the
following. The correspondence $T_p\times T_p$ on $\CC^2$ is given by the
map from $Y_0(p)\times Y_0(p)$ to $\CC^2\times\CC^2$ that sends a point
$(\phi,\psi)$ to $(s(\phi),s(\psi),t(\phi),t(\psi))$, where $s$ and $t$ stand
for source and target, respectively. Take the inverse image
of $C\times C$ in $Y_0(p)\times Y_0(p)$, and delete its zero-dimensional
part; that, together with its two maps to $C$, is~$T_{C,p}$.
We have to show that a suitable product
$T_{C,p_1}\cdots T_{C,p_r}$ with $r\geq1$ and the $p_i$ distinct induces a
non-trivial correspondence from $C_\CC$ to itself, because then we can apply
Theorem~\ref{thm6.1} to $C_\CC$ with $n=p_1\cdots p_r$. Let $S$ be the
finite set of irreducible components of~$C$. Then each $T_{C,p}$ induces
a correspondence $T_{S,p}$ on $S$ that is surjective in the sense that
both maps from $T_{S,p}$ to $S$ are surjective. Moreover, the Galois
group $G_\QQ$ acts transitively on $S$, and all $T_{S,p}$ are compatible
with this action. Let $x_0$ in $S$ correspond to~$C_\CC$. If there is some
$T_{S,p}$ such that $x_0$ is in $T_{S,p}x_0$, we can take $n=p$. So suppose
that for all $T_{S,p}$ we have $x_0\not\in T_{S,p}x_0$. Then we have
for all $T_{S,p}$ and all $x$ that $x\not\in T_{S,p}x$. One now easily
sees that there are $p_1,\ldots,p_r$ distinct with $1\leq r\leq|S|$ and
$x_0\in T_{S,p_1}\cdots T_{S,p_r}x_0$.
\begin{proof}
(Of Theorem~\ref{thm6.1}.) We take an integer $n$ as in the theorem we
are proving. Let $T_{C,n}$ be the correspondence on $C$ induced by
$T_n\times T_n$, in the sense explained above. (In fact, for everything that
follows we could also replace $T_{C,n}$ by one of its irreducible components,
but it is useful to see how to exploit all of~it.) We view
$T_{C,n}$ as a subset of $C\times C$. The image of $T_{C,n}$ under the
map $(\pr_1,\pr_1)$ from $C\times C$ to $\CC\times\CC$ is the image $T_n$
of $Y_0(n)$ in $\CC\times\CC$. Consider the commutative diagram:
\begin{eqn}\label{eqn6.2}
\renewcommand{\arraystretch}{1.5}
\begin{array}{ccc}
C & \to & \CC \\
\uparrow & & \uparrow \\
T_{C,n} & \to & T_n
\end{array}
\end{eqn}
in which the vertical maps are induced by the projections from $C\times C$
and $\CC\times\CC$ on the first factor.

\begin{lemma}\label{lemma6.3}
The map from $T_{C,n}$ to the fibred product $C\times_\CC T_n$ induced
by (\ref{eqn6.2}) is surjective.
\end{lemma}
\begin{proof}
By construction, all four maps in (\ref{eqn6.2}) are finite as morphisms of
(possibly reducible) algebraic curves. Therefore, the map from $T_{C,n}$
to $C\times_\CC T_n$ is also a finite morphism of algebraic curves.
Hence to show that it is surjective, it suffices to show that
$C\times_\CC T_n$ is irreducible, or, equivalently, that the tensor
product of the function fields of $C$ and $Y_0(n)$ over $\CC(j)$ is a field.
For this, it is enough to prove that the tensor product with $Y_0(n)$
replaced by $Y(n)$ is a field ($Y(n)$ is the modular curve parametrizing
elliptic curves with a symplectic basis of their $n$-torsion). The function
field of $Y(n)$ is Galois over $\CC(j)$ with Galois group
$\SL_2(\ZZ/n\ZZ)/\{\pm1\}$. The group $\SL_2(\ZZ/n\ZZ)$ is isomorphic to
the product of the $\SL_2(\FF_{p_i})$, $1\leq i\leq r$; one checks easily
that it has no non-trivial subgroup of index at most~$d_1$. This means that
the function fields of $C$ and $Y(n)$ are linearly disjoint.
\end{proof}
For reasons to become clear soon, we now first prove the following
lemma.
\begin{lemma}\label{lemma6.4}
The orbits in $C$ of $T_{C,n}$ are not discrete for the strong topology.
\end{lemma}
\begin{proof}
The morphism $\pr_1$ from $C$ to $\CC$ is proper, hence the image of a
closed subset of $C$ is closed in~$\CC$. In particular, the image of
the closure of any subset of $C$ is the closure of its image. Hence it
is enough to see that the images in $\CC$ of the orbits of $T_{C,n}$ are
not discrete. Let $x$ be in $C$, and let $y$ be its image in~$\CC$.
Lemma~\ref{lemma6.3} implies that $\pr_1 T_{C,n}x=T_ny$, hence we just
have to show that the orbits in $\CC$ of $T_n$ are not discrete. For this
we view $\CC$ as the quotient of the complex upper half plane $\HH$ by the
group $\SL_2(\ZZ)$ via the map $\pi\colon\tau\mapsto j(\CC/(\ZZ+\ZZ\tau))$.
Let $x$ be in $\CC$, and choose $\tau$ in~$\pi^{-1}x$. Then for all $a$ and
$b$ in $\ZZ$, $\pi(\tau+a)$ and $\pi(n^b\tau)$ are in the orbit of $x$
under~$T_n$. By composing these operations, we see that $\pi(n^b\tau+a)$
and $\pi(\tau+n^{-b}a)$ are in the orbit of~$x$. Taking $a$ non-zero and
$b$ big shows that the orbit is not discrete.
\end{proof}
We view $\CC\times\CC$ as the quotient of $\HH\times\HH$ by the
group $\Gamma:=\SL_2(\ZZ)\times\SL_2(\ZZ)$, via the map:
\begin{eqn}\label{eqn6.5}
\pi\colon \HH\times\HH \to \CC\times\CC, \quad
(\tau_1,\tau_2) \mapsto (j(\CC/(\ZZ+\ZZ\tau_1)),j(\CC/(\ZZ+\ZZ\tau_2))).
\end{eqn}
Let $X$ be an irreducible component of the analytic subvariety
$\pi^{-1}C$ of $\HH\times\HH$. The group $G:=\SL_2(\RR)\times\SL_2(\RR)$
acts transitively on~$\HH\times\HH$. We will study its subgroup $G_X$,
the stabilizer of~$X$. What we have to prove is that $G_X$ is the graph
of an inner automorphism of $\SL_2(\RR)$; this automorphism then tells us
for which $m$ our curve $C$ is the image
of~$Y_0(m)$. The decisive step in the proof of this is to see that $G_X$
is not discrete (if $C$ is an arbitrary curve in $\CC^2$, then $G_X$ is
typically discrete).

\begin{lemma}\label{lemma6.6}
The group $G_X$ is an analytic subgroup of $G$.
\end{lemma}
\begin{proof}
The action of $G$ on $\HH\times\HH$ is algebraic (it is given by fractional
linear transformations). The subgroup $G_X$ consists of exactly those
elements $g$ in $G$ that satisfy, for all $x$ in $X$, the two conditions
$gx\in X$ and $g^{-1}x\in X$. All these conditions are analytic.
\end{proof}
\begin{lemma}\label{lemma6.7}
The kernels of the two projections from $G_X$ to $\SL_2(\RR)$ are discrete.
\end{lemma}
\begin{proof}
This kernel $K$, say for the second projection, is the same as the stabilizer
of $X$ in the subgroup $\SL_2(\RR)\times\{1\}$ of~$G$. For all $\tau$ in
$\HH$, it stabilizes $X_\tau:=X\cap(\HH\times\{\tau\})$, which is discrete
since $d_2>0$; hence the connected component $K^o$ of $K$ stabilizes every
element of $X_\tau$. Now the stabilizer in $\SL_2(\RR)$ of the element
$i$ of $\HH$ is~$\SO_2(\RR)$. Because $d_1>0$, $K^o$ is contained in all
conjugates of $\SO_2(\RR)$, which is~$\{\pm1\}$.
\end{proof}

\begin{lemma}\label{lemma6.8}
The image in\/ $\SL_2(\ZZ)$ of\/ $\Gamma_X$, the stabilizer of\/ $X$ in\/
$\Gamma$, under the $i$th projection, has index at most~$d_i$.
\end{lemma}
\begin{proof}
We do the proof for $i=2$. We factor the map
$\pi\colon\HH\times\HH\to\CC\times\CC$ as follows:
\begin{subeqn}\label{eqn6.8.1}
\HH\times\HH \to \CC\times\HH \to \CC\times\CC .
\end{subeqn}
Let $Y$ be the image of $X$ in $\CC\times\HH$. Then $Y$ is an irreducible
component of the inverse image $Z$ of $C$ in~$\CC\times\HH$. The map from
$X$ to $C$ is the quotient for the action of $\Gamma_X$, hence the map
from $Y$ to $C$ is the quotient for the action of~$\pr_2\Gamma_X$.
It follows that $\pr_2\Gamma_X$ is the stabilizer in $\SL_2(\ZZ)$ of
$Y$ in $Z$, so the set $\SL_2(\ZZ)/\pr_2\Gamma_X$ is the set of irreducible
components of~$Z$. But $Z$ is also the fibred product of
$\pr_2\colon C\to\CC$ and $\HH\to\CC$, which implies that $Z$ has at most
$d_2$ irreducible components.
\end{proof}
Lemmas~\ref{lemma6.6}, \ref{lemma6.7} and \ref{lemma6.8} are in fact valid
for any curve $C$ in $\CC^2$ for which $d_1$ and $d_2$ are non-zero. The
next one crucially exploits that $C\subset (T_n\times T_n)C$.

\begin{lemma}\label{lemma6.9}
The topological group $G_X$ is not discrete.
\end{lemma}
\begin{proof}
The subgroup $G_X$ of $G$ is analytic, hence closed. It contains~$\Gamma_X$.
The inclusion $C\subset (T_n\times T_n)C$ implies that it contains some
less trivial elements as well. The correspondence $T_n$ on $\CC$ can be
described as follows. Take $z$ in~$\CC$; take its inverse image in~$\HH$;
apply the map $\tau\mapsto n\tau=({n\atop0}{0\atop1})\tau$ to it and take
its image in $\CC$; that is~$T_nz$. Another way to say this is: take
representatives in $\GL_2(\QQ)$ (there are $\psi(n)$ of them) $t_i$ for the
quotient set $\SL_2(\ZZ)({n\atop0}{0\atop1})\SL_2(\ZZ)/\SL_2(\ZZ)$; then for
$z$ in $\CC$ and $\tau$ in $\HH$ mapping to it, $T_nz$ is the image of the
sum of the~$t_i\tau$. It follows that for each $(i,j)$ such that $(t_i,t_j)X$
is contained in $\pi^{-1}C$ we get an element $g_{i,j}$ in $G_X$ of the form
$$
g_{i,j} = \gamma_{i,j,1}\cdot
\left(n^{-1/2}\left({n\atop0}{0\atop1}\right),
n^{-1/2}\left({n\atop0}{0\atop1}\right)\right)
\cdot\gamma_{i,j,2},
$$
with $\gamma_{i,j,1}$ and $\gamma_{i,j,2}$ in~$\Gamma$. For $c$ in $C$ and
$x$ in $X$ mapping to $c$, $T_{C,n}c$ is the image of the sum of
the~$g_{i,j}x$. Let $H$ be the subgroup
of $G_X$ generated by $\Gamma_X$ and these elements~$g_{i,j}$. We will prove
that $H$ is not discrete. Let $\Hbar$ be the closure of~$H$. We take an
element $x$ in~$X$. The map from $G$ to $\HH\times\HH$ sending $g$ to
$gx$ is proper, because the stabilizers of elements of $\HH\times\HH$ are
compact. Hence $\Hbar x$ is also the closure of~$Hx$. The subset $Hx$ of
$X$ is discrete if and only if its image in $C$ is discrete,
since $H$ contains $\Gamma_X$ and the map $X\to C$ is the quotient for the
action of~$\Gamma_X$. By construction, the image of $Hx$ in $C$ is the
orbit of $x$ for $T_{C,n}$, which, by Lemma~\ref{lemma6.4}, is not discrete.
This proves that $G_X$ is not discrete.
\end{proof}
We can now quickly finish the proof of Theorem~\ref{thm6.1}. Consider the
Lie algebra $\Lie(G_X)$, which by Lemma~\ref{lemma6.9} is non-zero.
Lemma~\ref{lemma6.7} tells us that the two projections $\pr_i\Lie(G_X)$
are non-zero. But $\pr_i\Lie(G_X)$ is normalized by $\pr_i\Gamma_X$, which
is Zariski dense in $\SL_2(\RR)$ by Lemma~\ref{lemma6.8}.
Since $\Lie(\SL_2(\RR))$ is simple, it follows that $\pr_i\Lie(G_X)$ is
equal to $\Lie(\SL_2(\RR))$ for both~$i$. So, since $\SL_2(\RR)$ is
connected, $G_X$ projects surjectively on both factors $\SL_2(\RR)$ of~$G$.
Now we apply what is called Goursat's lemma.
The kernel of $\pr_1\colon G_X\to\SL_2(\RR)$ is a normal subgroup of
$\SL_2(\RR)$, viewed as $\SL_2(\RR)\times\{1\}$. Since it is discrete and
contains $\{1,-1\}$, it is $\{1,-1\}$. The same holds for the other
projection, and $G_X$ is the inverse image in $G$ of the graph of an
analytic automorphism, $\sigma$ say, of~$\SL_2(\RR)/\{\pm1\}$. Every
such automorphism is inner. Since the $\pr_i\Gamma_X$ have finite index
in $\SL_2(\ZZ)$, it follows that $\sigma$ is induced from an inner
automorphism of the algebraic group~$\SL_{2,\QQ}$. The algebraic group of
automorphisms of $\SL_{2,\QQ}$ is~$\PGL_{2,\QQ}$. Since the map
$\GL_2(\QQ)\to\PGL_2(\QQ)$ is surjective (for example by Hilbert~90),
$\sigma$ is given by conjugation by some element $g$ in~$\GL_2(\QQ)$.
So $G_X$ is the set $\{(h,\pm ghg^{-1})\,|\, h\in \SL_2(\RR)\}$.
Let $x$ be an element of $X$, and write it as $x=(\tau,h\tau)$ with
$\tau$ in $\HH$ and $h$ in~$\SL_2(\RR)$. Since $G_Xx$ is in $X$, which
is of dimension two, the stabilizer of $x$ in $G_X$ has dimension at
least one. Let $H$ be the stabilizer of $\tau$ in the connected component
of identity $G_X^o$, for the action
of $G_X^o$ on the first factor $\HH$; then the stabilizer of $h\tau$ for the
action on the second factor is the conjugate $g^{-1}hHh^{-1}g$ of~$H$.
Since $H$ is of dimension one and connected (it is isomorphic to
$\SO_2(\RR)$) we must have $H=g^{-1}hHh^{-1}g$, i.e., $g^{-1}h$
normalizes~$H$. Since the normalizer of $\SO_2(\RR)$ in $\SL_2(\RR)$ is
just $\SO_2(\RR)$ itself, this means that $g^{-1}h$ is in $H$, or,
equivalently, that $h\tau=g\tau$. This means that
$X=\{(\tau,g\tau)\,|\,\tau\in\HH\}$. We may replace $g$ by multiples
$ag$ of it, with $a$ a non-zero rational number. So we can and do
suppose that $g\ZZ^2$ is contained in $\ZZ^2$ and that $\ZZ^2/g\ZZ^2$ is
cyclic, say of order~$m$. It is now clear that $C$ is~$Y_0(m)$.
\end{proof}

\section{Some remarks.}\label{section7}
\begin{remark}\label{rmk7.1}
Our proof of Theorem~\ref{thm1.1} shows in fact that, assuming GRH,
for each pair $(d_1,d_2)$ of positive integers there exists an effectively
computable number $B(d_1,d_2)$, such
that on every irreducible curve $C$ in $\CC^2$ of bi-degree $(d_1,d_2)$
that is defined over $\QQ$ and not a modular curve there are at most
$B(d_1,d_2)$ CM points.
(Note that under GRH, the statement that
$|\Pic(O_K)|/|\Pic(O_K)[2]|\to\infty$
is effective.)
\end{remark}

\begin{remark}\label{rmk7.2}
It is not true that all irreducible curves $C$ in $\CC^2$ with
$C\subset(T_n\times T_n)C$ for some $n>1$ are the image of some~$Y_0(m)$.
Here we construct some examples. Let $n>1$. Let $w_n$ be the Atkin-Lehner
involution of $Y_0(n)$: it sends an isogeny to its dual. The correspondence
$T_n$ on $\CC$ has the following description. For $z$ in $\CC$, take its
inverse image in $Y_0(n)$, take the image of that under $w_n$ and then the
image in~$\CC$. It follows that for an irreducible curve $C$ in $\CC^2$
such that at least one of the irreducible components of its inverse
image in $Y_0(n)\times Y_0(n)$ is stable under the involution $(w_n,w_n)$
we have $C\subset (T_n\times T_n)C$. Let $Z$ be the quotient of
$Y_0(n)\times Y_0(n)$ by that involution. Bertini's theorem,
see for example \cite[Theorem~6.3]{Jouanolou},
gives the existence of whole families of curves in $Z$ with irreducible
inverse image in $Y_0(n)\times Y_0(n)$. Take $C$ to be the image in $\CC^2$
of such an inverse image.
\end{remark}

\begin{remark}\label{rmk7.3}
The condition that $n$ be square free in Theorem~\ref{thm6.1} should not
be necessary; it is due to the laziness of the author.
\end{remark}

\begin{remark}\label{rmk7.4}
It is very tempting to try to generalize the methods of this article
to the general case of Oort's conjecture.
\end{remark}

\vspace{2\baselineskip}\noindent
{\bf Acknowledgements.} I would like to thank Rutger Noot for interesting
discussions on this subject that motivated me enough to work on it, and for
his remarks on previous versions of this article. I thank
Johan de Jong for his interest and the reference to~\cite{Serre1}. I am
very grateful to Etienne Fouvry for a letter in which he explains in detail
the results mentioned in Remark~\ref{rmk5.4}. Tim Dokshitzer pointed out
a gap in a previous version of this article. I want to thank Fabrice
Rouiller for helping me installing the necessary software on the computer
with which this article is written. Finally, I am very grateful for an
invitation to the Centre for Research in Mathematics at the Institut
d'Estudis Catalans in Barcelona, where I could compile my somewhat chaotic
and incomplete notes into this article.

\vfill
\noindent
Bas Edixhoven\\
IRMAR\\
Campus de Beaulieu\\
35042 Rennes cedex\\
France

\end{document}